\documentclass{article}
\usepackage{amssymb}

\usepackage{amsmath}

\input{tcilatex}

\begin{document}

\title{Special values of the Lommel functions and associated integrals}
\author{Bernard J. Laurenzi \\
Department of Chemistry\\
The State University of New York at Albany}
\date{December 11, 2018}
\maketitle

\begin{abstract}
Special values of the Lommel functions allow the calculation of Fresnel like
integrals. These closed form expressions along with their asymptotic values
are reported.
\end{abstract}

\section{Special values of the Lommel functions $s_{\protect\mu ,\protect\nu %
}(z)$ and $S_{\protect\mu ,\protect\nu }(z)$}

\ \ \ \ The Lommel functions $s_{\mu ,\nu }(z)$ and $S_{\mu ,\nu }(z)$ with
unrestricted $\mu ,\nu $ satisfy the relations \cite{Magnus1}

\begin{quote}
\begin{equation}
\frac{2\nu }{z}S_{\mu ,\nu }(z)=(\mu +\nu -1)S_{\mu -1,\nu -1}(z)-(\mu -\nu
-1)S_{\mu -1,\nu +1}(z),  \label{eq1}
\end{equation}

\begin{equation}
\lbrack (\mu +1)^{2}-\nu ^{2}]S_{\mu ,\nu }(z)+S_{\mu +2,\nu }(z)=z^{\mu +1},
\label{eq2}
\end{equation}%
\begin{equation}
\frac{dS_{\mu ,\nu }(z)}{dz}+\frac{\nu }{z}S_{\mu ,\nu }(z)=(\mu +\nu
-1)S_{\mu -1,\nu -1}(z),  \label{eq3}
\end{equation}%
and the symmetry property%
\begin{equation*}
S_{\mu ,-\nu }(z)=S_{\mu ,\nu }(z).
\end{equation*}

In the case where $\mu $ is an integer $k$ and $\nu =1/2$, the recurrence
relation (2) viewed as a second-order difference equation is 
\begin{equation}
(2k+1)(2k+3)S_{k,1/2}(z)+4S_{k+2,1/2}(z)=4z^{\,k+1},  \label{eq4}
\end{equation}%
and can be reduced to first-order equations in the cases of even or odd
values of $k.$ \ 

\bigskip
\end{quote}

\subsection{Case $k=2m$}

\begin{quote}
With $k=2m$ equation (4) can be written 
\begin{equation}
(4m+1)(4m+3)\,f_{m}(z)+4\,f_{m+1}(z)=4z^{\,2m+1},  \label{eq5}
\end{equation}

where $f_{m}(z)=S_{2m,1/2}(z)$ or $s_{2m,1/2}(z).$ It has the solution%
\begin{equation*}
f_{m}(z)=(-1)^{m}{\Gamma (2}m+1/2{)}\left[ \frac{f_{0}(z)}{\sqrt{\pi }}%
-z\sum_{j=0}^{m-1}\frac{(-z^{2})^{\,j}}{\Gamma (2j+5/2)}\right] ,
\end{equation*}%
\ the initial values of $\ f_{0}(z)$ being either $S_{0,1/2}(z)$ or $%
s_{0,1/2}(z).$ \ Using those relations the solutions to (5) become

\begin{subequations}
\label{0}
\begin{align}
S_{2m,1/2}(z)& =(-1)^{m}\Gamma (2m+1/2)[\tfrac{S_{0,1/2}(z)}{\sqrt{\pi }}%
-z\sum_{j=0}^{m-1}\tfrac{(-z^{2})^{\,j}}{\Gamma (2j+5/2)}],  \label{eq6a} \\
s_{2m,1/2}(z)& =(-1)^{m}\Gamma (2m+1/2)[\tfrac{s_{0,1/2}(z)}{\sqrt{\pi }}%
-z\sum_{j=0}^{m-1}\tfrac{(-z^{2})^{\,j}}{\Gamma (2j+5/2)}],  \label{eq6b}
\end{align}%
where $S_{0,1/2}(z),$ and $s_{0,1/2}(z)$ have been given by Magnus, \textit{%
et al} \cite{Magnus2} as 
\end{subequations}
\begin{eqnarray*}
\frac{1}{\sqrt{\pi }}S_{0,1/2}(z) &=&\sqrt{\frac{2}{z}}\left\{ \cos (z)[%
\frac{1}{2}-S(\chi )]-\sin (z)[\frac{1}{2}-C(\chi )]\right\} , \\
\frac{1}{\sqrt{\pi }}s_{0,1/2}(z) &=&\sqrt{\frac{2}{z}}\left\{ \sin
(z)C(\chi )-\cos (z)S(\chi )\right\} ,
\end{eqnarray*}%
with%
\begin{equation*}
\chi =\sqrt{\frac{2z}{\pi }},
\end{equation*}

and where $S$ and $C$ are the Fresnel sine and cosine integrals \cite{NBS}%
\begin{eqnarray*}
S(z) &=&\int_{0}^{z}\sin (\frac{1}{2}\pi t^{2})\,dt, \\
C(z) &=&\int_{0}^{z}\cos (\frac{1}{2}\pi t^{2})\,dt.
\end{eqnarray*}
\end{quote}

\subsection{Case $k=2m+1$}

\begin{quote}
\ \ \ \ In the case where $k=2m+1,$ the difference equation (2) becomes%
\begin{equation}
(4m+3)(4m+5)\,f_{m}(z)+4\,f_{m+1}(z)=4z^{\,2m+2},  \label{eq7}
\end{equation}

where $f_{m}(z)=S_{2m+1,1/2}(z)$ or $s_{2m+1,1/2}(z).$ Here the solution is 
\begin{equation*}
f_{m}(z)=(-1)^{m}\Gamma (2m+3/2)\left[ \frac{2\,f_{0}(z)}{\sqrt{\pi }}%
-z^{2}\sum_{j=0}^{m-1}\frac{(-z^{2})^{\,j}}{\Gamma (2j+7/2)}\right] .
\end{equation*}%
In the special case $\mu =-1$ and $\nu =1/2$ in (2) the functions\thinspace\ 
$f_{0}(z)$%
\begin{equation*}
f_{0}(z)=S_{1,1/2}(z)=1+\frac{1}{4}S_{-1,1/2}(z),
\end{equation*}%
or%
\begin{equation*}
f_{0}(z)=s_{1,1/2}(z)=1+\frac{1}{4}s_{-1,1/2}(z).
\end{equation*}%
We have%
\begin{eqnarray*}
S_{1,1/2}(z) &=&1+\sqrt{\frac{\pi }{2z}}\{\cos (z)[\frac{1}{2}-C(\chi
)]+\sin (z)[\frac{1}{2}-S(\chi )]\}, \\
s_{1,1/2}(z) &=&1-\sqrt{\frac{\pi }{2z}}\{\sin (z)S(\chi )+\cos (z)C(\chi
)\},
\end{eqnarray*}%
where values of $S_{-1,1/2}(z),$ $s_{-1,1/2}(z)$ have been given by Magnus, 
\textit{et al }as

\begin{eqnarray*}
\frac{\,1}{2\sqrt{\pi }}S_{-1,1/2}(z) &=&\sqrt{\frac{2}{z}}\left\{ \cos (z)[%
\frac{1}{2}-C(\chi )]+\sin (z)[\frac{1}{2}-S(\chi )]\right\} , \\
\frac{\,1}{2\sqrt{\pi }}s_{-1,1/2}(z) &=&-\sqrt{\frac{2}{z}}\left\{ \sin
(z)S(\chi )+\cos (z)C(\chi )\right\} .
\end{eqnarray*}
\end{quote}

Using those relations, the solutions to (7) become after resumming 
\begin{subequations}
\begin{eqnarray}
S_{2m+1,1/2}(z) &=&(-1)^{m}\Gamma (2m+3/2)[\frac{S_{-1,1/2}(z)}{2\sqrt{\pi }}%
+\sum_{j=0}^{m}\tfrac{(-z^{2})^{\,j}}{\Gamma (2j+3/2)}],  \label{eq8a} \\
s_{2m+1,1/2}(z) &=&(-1)^{m}\Gamma (2m+3/2)[\frac{s_{-1,1/2}(z)}{2\sqrt{\pi }}%
+\sum_{j=0}^{m}\tfrac{(-z^{2})^{\,j}}{\Gamma (2j+3/2)}].  \label{eq8b}
\end{eqnarray}

\subsection{Integrals with values containing the Fresnel functions}

The integrals 
\end{subequations}
\begin{eqnarray*}
&&\int_{0}^{1}z^{2k}\cos (\lambda z^{2})dz, \\
&&\int_{0}^{1}z^{2k}\sin (\lambda z^{2})dz,
\end{eqnarray*}%
which contain even powers of the variable, can be expressed in terms of the
Lommel functions. \ 

In the first instance, Maple gives 
\begin{eqnarray}
\int_{0}^{1}z^{2k}\cos (\lambda z^{2})dz &=&[1-\frac{s_{k+1,1/2(\lambda )}}{%
\lambda ^{k}}]\frac{\cos (\lambda )}{(2k+1)}  \label{eq9} \\
&&+[(2k-1)s_{k,3/2}+(2/\lambda )s_{k+1,1/2}]\frac{\sin (\lambda )}{2\lambda
^{k}(2k+1)}.  \notag
\end{eqnarray}%
Using (1) and (2) we get the simplified form%
\begin{equation}
\int_{0}^{1}z^{2k}\cos (\lambda z^{2})dz=\frac{1}{4\lambda ^{k}}[(2k-1)\cos
(\lambda )\,s_{k-1,1/2}(\lambda )+2\sin (\lambda )\,s_{k,1/2}(\lambda )].
\label{eq10}
\end{equation}%
The values of the integral in the case where $k=2m$ can be obtained from the
Lommel expressions above in (6b) and (8b) with $m$ replaced with $m-1$. \ We
get%
\begin{eqnarray*}
\int_{0}^{1}z^{4m}\cos (\lambda z^{2})dz &=&\frac{(-1)^{m}\Gamma (2m+1/2)}{%
2\lambda ^{2m}}\{\sqrt{\frac{2}{\lambda }}C(\chi ) \\
&&-\cos (\lambda )\sum_{j=0}^{m-1}\frac{(-\lambda ^{2})^{\,j}}{\Gamma
(2j+3/2)} \\
&&-\lambda \sin (\lambda )\sum_{j=0}^{m-1}\frac{(-\lambda ^{2})^{\,j}}{%
\Gamma (2j+5/2)}\}.
\end{eqnarray*}%
In the case $k=2m+1$ we have%
\begin{eqnarray*}
\int_{0}^{1}z^{4m+2}\cos (\lambda z^{2})dz &=&\frac{(-1)^{m+1}\Gamma (2m+3/2)%
}{2\lambda ^{2m+1}}\{\sqrt{\frac{2}{\lambda }}S(\chi ) \\
&&+\lambda \cos (\lambda )\sum_{j=0}^{m-1}\frac{(-\lambda ^{2})^{\,j}}{%
\Gamma (2j+5/2)} \\
&&-\sin (\lambda )\sum_{j=0}^{m}\frac{(-\lambda ^{2})^{\,j}}{\Gamma (2j+3/2)}%
\}.
\end{eqnarray*}%
With these results we see that all cases of cosine integrals with even
powers of the variable have been obtained in terms containing the Fresnel $S$
function. \ 

For the corresponding sine integrals, integration by parts gives%
\begin{equation*}
\int_{0}^{1}z^{2k}\sin (\lambda z^{2})dz=\frac{\sin (\lambda )}{2k+1}-\frac{%
2\lambda }{2k+1}\int_{0}^{1}z^{2(k+1)}\cos (\lambda z^{2})dz.
\end{equation*}%
Using (10) we have%
\begin{equation*}
\int_{0}^{1}z^{2k}\sin (\lambda z^{2})dz=\frac{1}{4\lambda ^{k}}[(2k-1)\sin
(\lambda )s_{k-1,1/2}(\lambda )-2\cos (\lambda )s_{k,1/2}(\lambda )].
\end{equation*}%
With $k=2m$ we get%
\begin{equation*}
\int_{0}^{1}z^{4m}\sin (\lambda z^{2})dz=\frac{1}{2\lambda ^{2m}}[(4m-1)\sin
(\lambda )s_{2m-1,1/2}(\lambda )-2\cos (\lambda )s_{2m,1/2}(\lambda )],
\end{equation*}%
which then becomes%
\begin{eqnarray*}
\int_{0}^{1}z^{4m}\sin (\lambda z^{2})dz &=&\frac{(-1)^{m}\Gamma (2m+1/2)}{%
2\lambda ^{2m}}\{\sqrt{\frac{2}{\lambda }}S(\chi ) \\
&&-\sin (\lambda )\sum_{j=0}^{m-1}\frac{(-1)^{\,j}\lambda ^{2j}}{\Gamma
(2j+3/2)} \\
&&+\lambda \cos (\lambda )\sum_{j=0}^{m-1}\frac{(-1)^{\,j}\lambda ^{2j}}{%
\Gamma (2j+5/2)}\}.
\end{eqnarray*}%
For $k=2m+1$ the sine integrals are given by Maple as%
\begin{equation*}
\int_{0}^{1}z^{4m+2}\sin (\lambda z^{2})dz=\frac{1}{4\lambda ^{2m+1}}[\sin
(\lambda )(4m+1)s_{2m,1/2}(\lambda )-2\cos (\lambda )s_{2m+1,1/2}(\lambda )].
\end{equation*}%
Using the values of Lommel functions given above we have%
\begin{eqnarray*}
\int_{0}^{1}z^{4m+2}\sin (\lambda z^{2})dz &=&\frac{(-1)^{m}\Gamma (2m+3/2)}{%
2\lambda ^{2m+1}}\{\sqrt{\frac{2}{\lambda }}C(\chi ) \\
&&-\lambda \sin (\lambda )\sum_{j=0}^{m-1}\frac{(-1)^{\,j}\lambda ^{2j}}{%
\Gamma (2j+5/2)} \\
&&-\cos (\lambda )\sum_{j=0}^{m}\frac{(-1)^{\,j}\lambda ^{2j}}{\Gamma
(2j+3/2)}\}.
\end{eqnarray*}%
With these results we see that all values of the sine integrals which
contain \textit{even} powers of $z$ requires the presence of the Fresnel
integrals $C(\chi )$.

\section{Integrals with values containing elementary functions}

We next consider integrals which contain \textit{odd} powers of $z$ i.e.%
\begin{eqnarray*}
&&\int_{0}^{1}z^{2k+1}\cos (\lambda z^{2})dz, \\
&&\int_{0}^{1}z^{2k+1}\sin (\lambda z^{2})dz.
\end{eqnarray*}%
In the first instance Maple gives%
\begin{eqnarray*}
\int_{0}^{1}z^{2k+1}\cos (\lambda z^{2})dz &=&\frac{1}{2(k+1)}\{(1-\frac{%
s_{k+3/2,1/2}(\lambda )}{\lambda ^{k+1/2}})\cos (\lambda ) \\
&&+(\,k\,s_{k+1/2,3/2}(\lambda ))+\frac{s_{k+3/2,1/2}(\lambda )}{\lambda })%
\frac{\sin (\lambda )}{\lambda ^{k+1/2}}\},
\end{eqnarray*}%
which reduces to%
\begin{equation}
\int_{0}^{1}z^{2k+1}\cos (\lambda z^{2})dz=\frac{1}{2\lambda ^{k+1/2}}[k\cos
(\lambda )s_{k-1/2,1/2}(\lambda )+\sin (\lambda )s_{k+1/2,1/2}(\lambda )]
\label{eq11}
\end{equation}%
using (1) and (2). \ Where $k=2m$ we have%
\begin{equation}
\int_{0}^{1}z^{4m+1}\cos (\lambda z^{2})dz=\frac{1}{2\lambda ^{2m+1/2}}%
[2m\cos (\lambda )s_{2m-1/2,1/2}(\lambda )+\sin (\lambda
)s_{2m+1/2,1/2}(\lambda )].  \label{eq12}
\end{equation}%
We see that this expression requires the Lommel functions $%
s_{2m-1/2,1/2}(\lambda )$ and $s_{2m+1/2,1/2}(\lambda )$. \ In the first
case we have from (2) and the initial condition $s_{3/2,1/2}(\lambda )=\sqrt{%
\lambda }\left[ 1-\sin (\lambda )/\lambda \right] $, the Lommel function $%
s_{2m-1/2,1/2}(\lambda )$ as 
\begin{equation*}
s_{2m-1/2,1/2}(\lambda )=\frac{(-1)^{m}(2m-1)!}{\sqrt{\lambda }}\left[ \sin
(\lambda )-\sum_{j=0}^{m-1}\frac{(-1)^{\,j}\lambda ^{2j+1}}{(2j+1)!}\right] .
\end{equation*}%
\ The function $s_{2m+1,1/2}(\lambda )$ then can be obtained from the
differential-difference equation in (3). \ That is to say%
\begin{equation}
\frac{d\,s_{2m+1/2,1/2}(\lambda )}{d\lambda }+\frac{1}{2\lambda }%
s_{2m+1/2,1/2}(\lambda )=2m\,s_{2m-1/2,1/2}(\lambda ).  \label{eq13}
\end{equation}%
We get with the initial condition $s_{2m+1/2,1/2}(0)=0,$ the solution to
(13) as%
\begin{equation}
s_{2m+1/2,1/2}(\lambda )=\frac{2m}{\sqrt{\lambda }}\int_{0}^{\lambda }\sqrt{z%
}s_{2m-1/2,1/2}(z)\,dz,  \label{eq14}
\end{equation}%
which immediately gives 
\begin{equation*}
s_{2m+1/2,1/2}(\lambda )=\frac{(-1)^{m+1}(2m)!}{\sqrt{\lambda }}\left[ \cos
(\lambda )-\sum_{j=0}^{m}\frac{(-1)^{\,j}\lambda ^{2j}}{(2j)!}\right] .
\end{equation*}%
We have obtained in these cases closed forms for the Lommel functions which
contain only elementary functions. As a result, we get%
\begin{eqnarray*}
\int_{0}^{1}z^{4m+1}\cos (\lambda z^{2})dz &=&\frac{(-1)^{m}(2m)!}{2\lambda
^{2m+1}}\{\sin (\lambda )\sum_{j=0}^{m}\frac{(-1)^{\,j}\lambda ^{2j}}{(2j)!}
\\
&&-\lambda \cos (\lambda )\sum_{j=0}^{m-1}\frac{(-1)^{\,j}\lambda ^{2j}}{%
(2j+1)!}\}.
\end{eqnarray*}%
The cosine integral containing powers $4m+3$ is given by%
\begin{equation*}
\int_{0}^{1}z^{4m+3}\cos (\lambda z^{2})dz=\frac{1}{2\lambda ^{2m+3/2}}%
[(2m+1)\cos (\lambda )\,s_{2m+1/2,1/2}(\lambda )+\sin (\lambda
)s_{2m+3/2,1/2}(\lambda )].
\end{equation*}%
In the latter expression the quantity $s_{2m+3/2,1/2}(\lambda )$ can be
obtained from (3) and we have 
\begin{eqnarray*}
\int_{0}^{1}z^{4m+3}\cos (\lambda z^{2})dz &=&\frac{(-1)^{m+1}(2m+1)!}{%
2\lambda ^{2m+2}}\{1-\cos (\lambda )\sum_{j=0}^{m}\frac{(-1)^{\,j}\lambda
^{2j}}{(2j)!} \\
&&-\lambda \sin (\lambda )\sum_{j=0}^{m}\frac{(-1)^{\,j}\lambda ^{2j}}{%
(2j+1)!}\}.
\end{eqnarray*}%
An expression for the sine integrals with odd powers i.e.%
\begin{equation*}
\int_{0}^{1}z^{2k+1}\sin (\lambda z^{2})dz=\frac{1}{2\lambda ^{k+1/2}}[k\sin
(\lambda )\,\,s_{k-1/2,1/2}(\lambda )-\cos (\lambda )s_{k+1/2,1/2}(\lambda
)],
\end{equation*}%
has been obtained using integration by parts of the corresponding $\cos
(\lambda z^{2})$ integral together with the latter's integrated form. \ Then
in the case where $\ k=2m$ we have 
\begin{eqnarray*}
\int_{0}^{1}z^{4m+1}\sin (\lambda z^{2})dz &=&\frac{(-1)^{m}(2m)!}{2\lambda
^{2m+1}}\{1-\cos (\lambda )\sum_{j=0}^{m}\frac{(-1)^{\,j}\lambda ^{2j}}{(2j)!%
} \\
&&-\lambda \sin (\lambda )\sum_{j=0}^{m-1}\frac{(-1)^{\,j}\lambda ^{2j}}{%
(2j+1)!}\}.
\end{eqnarray*}%
Where $k=2m+1$ we have%
\begin{equation*}
\int_{0}^{1}z^{4m+3}\sin (\lambda z^{2})dz=\frac{1}{2\lambda ^{2m+3/2}}%
[(2m+1)\sin (\lambda )\,\,s_{2m+1/2,1/2}(\lambda )-\cos (\lambda
)s_{2m+3/2,1/2}(\lambda ).
\end{equation*}%
In the latter expression the quantity $s_{2m+3/2,1/2}(\lambda )$ can also be
obtained from (3) and the integral being sought has the value

\begin{eqnarray*}
\int_{0}^{1}z^{4m+3}\sin (\lambda z^{2})dz &=&\frac{(-1)^{m}(2m+1)!}{%
2\lambda ^{2m+2}}\{-\lambda \cos (\lambda )\sum_{j=0}^{m}\frac{%
(-1)^{\,j}\lambda ^{2j}}{(2j+1)!} \\
&&+\sin (\lambda )\sum_{j=0}^{m}\frac{(-1)^{\,j}\lambda ^{2j}}{(2j)!}\}.
\end{eqnarray*}%
\ 

\section{Asymptotic forms for the integrals containing $S(\protect\chi )$
and $C(\protect\chi )$}

It is useful to provide values of the integrals evaluated above in section 2
where the parameter $\lambda $ is large. \ Initially we consider the case of
the integrals with even powers of the integration variable $z$ i.e.%
\begin{eqnarray*}
&&\int_{0}^{1}z^{4m}\cos (\lambda z^{2})dz, \\
&&\int_{0}^{1}z^{4m}\sin (\lambda z^{2})dz, \\
&&\int_{0}^{1}z^{4m+2}\cos (\lambda z^{2})dz, \\
&&\int_{0}^{1}z^{4m+2}\sin (\lambda z^{2})dz,
\end{eqnarray*}%
which we have seen contain the Fresnel integrals. The asymptotic forms for
the functions $S\left( \chi \right) $ and $C(\chi )$ can be obtained from
the expressions \cite{NIST} 
\begin{eqnarray*}
S\left( \sqrt{\frac{2\lambda }{\pi }}\right) &=&\frac{1}{2}-f\,(\sqrt{\frac{%
2\lambda }{\pi }})\cos (\lambda )-g(\sqrt{\frac{2\lambda }{\pi }})\sin
(\lambda ), \\
C\left( \sqrt{\frac{2\lambda }{\pi }}\right) &=&\frac{1}{2}+f\,(\sqrt{\frac{%
2\lambda }{\pi }})\sin (\lambda )-g(\sqrt{\frac{2\lambda }{\pi }})\cos
(\lambda ),
\end{eqnarray*}%
where $f(z)$ and $g(z)$ are the Fresnel auxiliary functions. \ Their
asymptotic forms with (cut of\thinspace f $N\geq 1$) are given by%
\begin{eqnarray*}
f\,(\sqrt{\frac{2\lambda }{\pi }}) &\thicksim &\frac{1}{\sqrt{2}\,\pi
\lambda ^{1/2}}\sum_{j=0}^{N-1}\frac{(-1)^{j}}{\lambda ^{2j}}\Gamma (2j+1/2),
\\
g\,(\sqrt{\frac{2\lambda }{\pi }}) &\thicksim &\frac{1}{\sqrt{2}\,\pi
\lambda ^{3/2}}\sum_{j=0}^{N-1}\frac{(-1)^{j}}{\lambda ^{2j}}\Gamma (2j+3/2).
\end{eqnarray*}%
We get for the integrals in question, having used the relation%
\begin{equation*}
\frac{1}{\Gamma (-z)}=-\frac{\Gamma (z+1)}{\pi }\sin (\pi z),
\end{equation*}%
in the sine and cosine sums, the expressions%
\begin{eqnarray*}
\int_{0}^{1}z^{4m}\cos (\lambda z^{2})dz &\thicksim &\frac{\Gamma (2m+1/2)}{2%
}\{\frac{(-1)^{m}}{\sqrt{2}\lambda ^{2m+1/2}} \\
&&+\frac{\cos (\lambda )}{\pi }\sum_{j=1}^{m+N}\frac{(-1)^{\,j}}{\lambda
^{2j}}\Gamma (2j-2m-1/2) \\
&&+\frac{\sin (\lambda )}{\pi }\sum_{j=1}^{m+N}\frac{(-1)^{\,j+1}}{\lambda
^{2j-1}}\Gamma (2j-2m-3/2)\},
\end{eqnarray*}%
\begin{eqnarray*}
\int_{0}^{1}z^{4m}\sin (\lambda z^{2})dz &\thicksim &\frac{\Gamma (2m+1/2)}{2%
}\{\frac{(-1)^{m}}{\sqrt{2}\lambda ^{2m+1/2}} \\
&&+\frac{\cos (\lambda )}{\pi }\sum_{j=1}^{m+N}\frac{(-1)^{\,j}}{\lambda
^{2j-1}}\Gamma (2j-2m-3/2) \\
&&+\frac{\sin (\lambda )}{\pi }\sum_{j=1}^{m+N}\frac{(-1)^{\,j}}{\lambda
^{2j}}\Gamma (2j-2m-1/2)\},
\end{eqnarray*}%
\begin{eqnarray*}
\int_{0}^{1}z^{4m+2}\cos (\lambda z^{2})dz &\thicksim &\frac{\Gamma (2m+3/2)%
}{2}\{\frac{(-1)^{m+1}}{\sqrt{2}\lambda ^{2m+3/2}} \\
&&+\frac{\cos (\lambda )}{\pi }\sum_{j=1}^{m+N}\frac{(-1)^{\,j+1}}{\lambda
^{2j}}\Gamma (2j-2m-3/2) \\
&&+\frac{\sin (\lambda )}{\pi }\sum_{j=1}^{m+N+1}\frac{(-1)^{\,j}}{\lambda
^{2j-1}}\Gamma (2j-2m-5/2)\},
\end{eqnarray*}%
\begin{eqnarray*}
\int_{0}^{1}z^{4m+2}\sin (\lambda z^{2})dz &\thicksim &\frac{\Gamma (2m+3/2)%
}{2}\{\frac{(-1)^{m}}{\sqrt{2}\lambda ^{2m+3/2}} \\
&&+\frac{\cos (\lambda )}{\pi }\sum_{j=1}^{m+N}\frac{(-1)^{\,j+1}}{\lambda
^{2j+1}}\Gamma (2j-2m-1/2) \\
&&+\frac{\sin (\lambda )}{\pi }\sum_{j=1}^{m+N}\frac{(-1)^{\,j+1}}{\lambda
^{2j}}\Gamma (2j-2m-3/2)\}.
\end{eqnarray*}%
We see that the terms in these asymptotic expressions contain the sine and
cosine sums along with the additional and significant contributions from the
Fresnel integrals.

\subsection{Asymptotic values for Fresnel related integrals}

\bigskip

In extensions of the Schwinger-Englert semi-classical theory \cite{englert}\
of atomic structure, asymptotic forms for the integrals%
\begin{eqnarray*}
&&\int_{0}^{1}\frac{\sin (\lambda z^{2})}{(1+a\,z^{2})}dz, \\
&&\int_{0}^{1}\frac{\cos (\lambda z^{2})}{(1+a\,z^{2})}dz,
\end{eqnarray*}%
arise with $\lambda \rightarrow \infty $ and $0<a<1.$ \ Here we give
estimates for those forms. \ Expanding the denominators of the integrals
above we have 
\begin{equation*}
\int_{0}^{1}\frac{\sin (\lambda z^{2})}{(1+a\,z^{2})}dz=\sum_{k=0}^{\infty
}a^{2k}\int_{0}^{1}z^{4k}\sin (\lambda z^{2})-a\sum_{k=0}^{\infty
}a^{2k}\int_{0}^{1}z^{4k+2}\sin (\lambda z^{2}),
\end{equation*}%
using the results obtained above we get (where $z!!$ is the double factorial
function and $N$ is the order of truncation)%
\begin{eqnarray*}
\int_{0}^{1}\frac{\sin (\lambda z^{2})}{(1+a\,z^{2})}dz &\thicksim &\frac{1}{%
2}\sqrt{\frac{\pi }{2\lambda }}\sum_{k=0}^{\infty }(-1)^{k}\left( \frac{a}{%
2\lambda }\right) ^{2k}[(4k-1)!!-\frac{a}{2\lambda }(4k+1)!!] \\
&&+\frac{\cos (\lambda )}{2}\{\sum_{k=0}^{\infty
}a^{2k}(4k-1)!!\sum_{j=1}^{k+N}\frac{(-1)^{j}(4j-4k-5)!!}{\left( 2\lambda
\right) ^{2j-1}} \\
&&+\frac{a}{\lambda }\sum_{k=0}^{\infty }a^{2k}(4k+1)!!\sum_{j=1}^{k+N}\frac{%
(-1)^{j}(4j-4k-3)!!}{\left( 2\lambda \right) ^{2j}}\} \\
&&+\frac{\sin (\lambda )}{2}\{\sum_{k=0}^{\infty
}a^{2k}(4k-1)!!\sum_{j=1}^{k+N}\frac{(-1)^{j}(4j-4k-3)!!}{\left( 2\lambda
\right) ^{2j}} \\
&&+\frac{a}{\lambda }\sum_{k=0}^{\infty }a^{2k}(4k+1)!!\sum_{j=1}^{k+N}\frac{%
(-1)^{j}(4j-4k-5)!!}{\left( 2\lambda \right) ^{2j-1}}\}.
\end{eqnarray*}%
In lowest order in $a$ and $1/\lambda $ we have%
\begin{equation}
\int_{0}^{1}\frac{\sin (\lambda z^{2})}{(1+a\,z^{2})}dz\thicksim \frac{1}{4}%
\sqrt{\frac{2\pi }{\lambda }}\left( 1-\frac{a}{2\lambda }\right) -\frac{\cos
(\lambda )}{2\lambda }(1-a)-\frac{\sin (\lambda )}{\left( 2\lambda \right)
^{2}}(1+a)+\cdots ,  \label{eq15}
\end{equation}%
and in the case of the cosine integral we have%
\begin{eqnarray*}
\int_{0}^{1}\frac{\cos (\lambda z^{2})}{(1+a\,z^{2})}dz &\thicksim &\frac{1}{%
2}\sqrt{\frac{\pi }{2\,\lambda }}\sum_{k=0}^{\infty }(-1)^{k}\left( \frac{a}{%
2\lambda }\right) ^{2k}[(4k-1)!!+\frac{a}{2\lambda }(4k+1)!!] \\
&&+\cos (\lambda )\{\sum_{k=0}^{\infty }a^{2k}(4k-1)!!\sum_{j=1}^{k+N}\frac{%
(-1)^{j}(4j-4k-3)!!}{\left( 2\lambda \right) ^{2j}} \\
&&+\sum_{k=0}^{\infty }a^{2k+1}(4k+1)!!\sum_{j=1}^{k+N}\frac{%
(-1)^{j}(4j-4k-5)!!}{\left( 2\lambda \right) ^{2j}}\} \\
&&-\sin (\lambda )\{\sum_{k=0}^{\infty }a^{2k}(4k-1)!!\sum_{j=1}^{k+N}\frac{%
(-1)^{j}(4j-4k-5)!!}{\left( 2\lambda \right) ^{2j-1}} \\
&&+\sum_{k=0}^{\infty }a^{2k+1}(4k+1)!!\sum_{j=1}^{k+N+1}\frac{%
(-1)^{j}(4j-4k-7)!!}{\left( 2\lambda \right) ^{2j-1}}\}.
\end{eqnarray*}%
In lowest order in $a$ and $1/\lambda $ we have%
\begin{equation}
\int_{0}^{1}\frac{\cos (\lambda z^{2})}{(1+a\,z^{2})}dz\thicksim \frac{1}{4}%
\sqrt{\frac{2\pi }{\lambda }}\left( 1+\frac{a}{2\lambda }\right) -\frac{\cos
(\lambda )}{\left( 2\lambda \right) ^{2}}(1+a)+\frac{\sin (\lambda )}{%
2\lambda }(1-a)+\cdots  \label{eq16}
\end{equation}%
It is interesting to note that the integrals with infinite range i.e.%
\begin{equation}
\int_{0}^{\infty }\frac{\sin (\lambda z^{2})}{(1+a\,z^{2})}=-\frac{\pi }{2%
\sqrt{a}}\sin (\lambda /a)+\frac{\pi }{2\sqrt{a}}\left\{ 
\begin{array}{c}
C(\sqrt{\frac{2\lambda }{\pi a}})\left[ \sin (\frac{\lambda }{a})+\cos (%
\frac{\lambda }{a})\right] \\ 
+S(\sqrt{\frac{2\lambda }{\pi a}})\left[ \sin (\frac{\lambda }{a})-\cos (%
\frac{\lambda }{a})\right]%
\end{array}%
\right\} ,  \label{eq17}
\end{equation}%
and%
\begin{equation}
\int_{0}^{\infty }\frac{\cos (\lambda z^{2})}{(1+a\,z^{2})}=\frac{\pi }{2%
\sqrt{a}}\cos (\lambda /a)-\frac{\pi }{2\sqrt{a}}\left\{ 
\begin{array}{c}
C(\sqrt{\frac{2\lambda }{\pi a}})\left[ \cos (\frac{\lambda }{a})-\sin (%
\frac{\lambda }{a})\right] \\ 
+S(\sqrt{\frac{2\lambda }{\pi a}})\left[ \cos (\frac{\lambda }{a}))+\sin (%
\frac{\lambda }{a})\right]%
\end{array}%
\right\} ,  \label{eq18}
\end{equation}%
which have the exact values shown above, posses\ the same forms for large $%
\lambda $ as the integrals (albeit without oscillations) with finite range.
\ That is to say 
\begin{equation}
\int_{0}^{\infty }\frac{\sin (\lambda z^{2})}{(1+a\,z^{2})}\thicksim \frac{1%
}{4}\sqrt{\frac{2\pi }{\lambda }}\left( 1-\frac{a}{2\lambda }\right) +\cdots
,  \label{eq19}
\end{equation}%
\begin{equation}
\int_{0}^{\infty }\frac{\cos (\lambda z^{2})}{(1+a\,z^{2})}dz\thicksim \frac{%
1}{4}\sqrt{\frac{2\pi }{\lambda }}\left( 1+\frac{a}{2\lambda }\right) +\cdots
\label{eq20}
\end{equation}%
\ \ \ 

It is also worth noting that asymptotic values of the integrals which
contain higher powers of the trigonometric function described above and with
higher powers of the denominators i.e. 
\begin{equation*}
\int_{0}^{1}\frac{\sin ^{n}(\lambda z^{2})}{(1+a\,z^{2})^{\,\nu }},\hspace{%
0.25in}\int_{0}^{1}\frac{\cos ^{n}(\lambda z^{2})}{(1+a\,z^{2})^{\,\nu }},
\end{equation*}%
and%
\begin{equation*}
\int_{0}^{\infty }\frac{\sin ^{n}(\lambda z^{2})}{(1+a\,z^{2})^{\,\nu }},%
\hspace{0.25in}\int_{0}^{\infty }\frac{\cos ^{n}(\lambda z^{2})}{%
(1+a\,z^{2})^{\,\nu }},
\end{equation*}%
can be expressed in terms of the forms given in (15,16) and (19,20). \ By
way of example, writing 
\begin{equation*}
I_{\nu }^{\,(\eta )}(a,\lambda )=\int_{0}^{\infty }\frac{\cos ^{\,\eta
}(\lambda z^{2})}{(1+a\,z^{2})^{\,\nu }},
\end{equation*}%
(noting the scaling property)%
\begin{equation*}
I_{\nu }^{\,(\eta )}(a,\lambda )=\frac{1}{\sqrt{a}}I_{\nu }^{\,(\eta
)}(1,\lambda /a),
\end{equation*}%
we get the differential-difference equation

\begin{equation}
I_{\,\nu +1\,}^{\,(\eta )}(a,\lambda )=I_{\nu }^{\,(\eta )}(a,\lambda
)+\left( \frac{a}{\nu }\right) \frac{d\,I_{\nu }^{\,(\eta )}(a,\lambda )}{da}%
.  \label{eq21}
\end{equation}%
In the case where $\eta =1$ and $\nu =1/2$ we have%
\begin{eqnarray*}
I_{1/2}^{\,(1)}(a,\lambda ) &=&\int_{0}^{\infty }\frac{\cos (\lambda z^{2})}{%
(1+a\,z^{2})^{\,1/2}} \\
&=&\frac{\pi }{4\sqrt{a}}\left\{ \sin (\frac{\lambda }{2a})\,J_{0}(\frac{%
\lambda }{2a})-\cos (\frac{\lambda }{2a})Y_{0}(\frac{\lambda }{2a})\right\} ,
\end{eqnarray*}%
a closed form expression where $J_{0}(z)$ and $Y_{0}(z)$ are Bessel
functions of the first and second kind. \ The value for this integral for
large $\frac{\lambda }{2a}$ using Hankel's \cite{Hankel}\ \ asymptotic
expressions for the Bessel functions is%
\begin{equation*}
I_{1/2}^{\,(1)}(a,\lambda )\thicksim \frac{1}{4}\sqrt{\frac{2\pi }{\lambda }}%
\left\{ 1+\frac{1}{4}\left( \frac{a}{\lambda }\right) -\frac{9}{32}\left( 
\frac{a}{\lambda }\right) ^{2}-\frac{75}{128}\left( \frac{a}{\lambda }%
\right) ^{3}+\cdots \right\} .
\end{equation*}%
In the case of the integral $I_{1}^{\,(2)}(a,\lambda )$, its exact value is 
\begin{eqnarray*}
I_{1}^{\,(2)}(a,\lambda ) &=&\int_{0}^{\infty }\frac{\cos ^{2}(\lambda z^{2})%
}{(1+a\,z^{2})} \\
&=&\frac{\pi }{2\sqrt{a}}\cos ^{2}(\lambda /a)-\frac{\pi }{4\sqrt{a}}\left\{ 
\begin{array}{c}
C(\sqrt{\frac{4\lambda }{\pi a}})\left[ \cos (\frac{2\lambda }{a})-\sin (%
\frac{2\lambda }{a})\right] \\ 
+S(\sqrt{\frac{4\lambda }{\pi a}})\left[ \cos (\frac{2\lambda }{a}))+\sin (%
\frac{2\lambda }{a})\right]%
\end{array}%
\right\} ,
\end{eqnarray*}%
or rewriting it in terms of the Fresnel auxiliary functions \cite{Aux}$\
f\,(\chi ^{\prime })$ and $g(\chi ^{\prime })$ referred to above (here with $%
\chi ^{\prime }=2\sqrt{\lambda /\pi a}$ ) we get 
\begin{equation*}
\int_{0}^{\infty }\frac{\cos ^{2}(\lambda z^{2})}{(1+a\,z^{2})}=\frac{\pi }{4%
\sqrt{a}}\left\{ 1+f\,(\chi ^{\prime })+g(\chi ^{\prime })\right\} .
\end{equation*}%
The asymptotic values of $\ f\,(\chi ^{\prime })$ and $\ g(\chi ^{\prime })$
for large $\lambda $ \cite{asymptotic} give the result%
\begin{eqnarray*}
\int_{0}^{\infty }\frac{\cos ^{2}(\lambda z^{2})}{(1+a\,z^{2})} &\thicksim &%
\frac{\pi }{4\sqrt{a}}+\frac{1}{8}\sqrt{\frac{\pi }{\lambda }}\left\{ 1-%
\frac{3}{16}\left( \frac{a}{\lambda }\right) ^{2}+\cdots \right\} \\
&&+\frac{1}{32}\sqrt{\frac{\pi }{\lambda }}\left\{ \frac{a}{\lambda }-\frac{%
15}{16}\left( \frac{a}{\lambda }\right) ^{3}+\cdots \right\} .
\end{eqnarray*}%
The first of the higher order integrals i.e. $\ I_{2}^{\,(2)}$ has the
closed form expression 
\begin{eqnarray*}
I_{2}^{\,(2)} &=&\int_{0}^{\infty }\frac{\cos ^{2}(\lambda z^{2})}{%
(1+a\,z^{2})^{2}} \\
&=&\frac{\pi }{8\sqrt{a}}\left\{ 
\begin{array}{c}
2\cos ^{2}(\frac{\lambda }{a})+4(\frac{\lambda }{a})\sin (\frac{2\lambda }{a}%
)+\sqrt{\frac{4\lambda }{\pi a}} \\ 
+C(\sqrt{\frac{4\lambda }{\pi a}})[\sin (\frac{2\lambda }{a})(1-4\frac{%
\lambda }{a})-\cos (\frac{2\lambda }{a})(1+4\frac{\lambda }{a})] \\ 
-S(\sqrt{\frac{4\lambda }{\pi a}})[\sin (\frac{2\lambda }{a})(1+4\frac{%
\lambda }{a})+\cos (\frac{\lambda }{a})(1-4\frac{2\lambda }{a}]%
\end{array}%
\right\} .
\end{eqnarray*}

The higher order integrals i.e. those containing powers of the cosine
function i.e. $\cos ^{\,\eta }(\lambda x^{2})$ are expressible in terms of
the integrals in (18). \ Using the trigonometric relations 
\begin{equation*}
\cos ^{\eta }(\lambda x^{2})=\frac{1}{2^{\eta }}\sum_{j=0}^{\eta }\binom{%
^{\eta }}{j}\,\cos ([\eta -2j]\lambda x^{2}),
\end{equation*}%
we have 
\begin{eqnarray*}
I_{\nu }^{\,(2\eta )}(a,\lambda ) &=&\sqrt{\frac{\pi }{a}}\frac{\Gamma
(2\eta -1/2)}{2^{\,2\eta }\,\eta !\,(\eta -1)!}+\frac{1}{2^{2\eta -1}}%
\sum_{j=1}^{\eta }\binom{2\,\eta }{\eta -j}\,I_{\nu }^{\,(1)}(a,2j\,\lambda
), \\
I_{\nu }^{\,(2\eta +1)}(a,\lambda ) &=&\frac{1}{2^{2\eta }}\sum_{j=0}^{\eta }%
\binom{2\eta +1}{\eta -j}\,I_{\nu }^{\,(1)}(a,(2j+1)\,\lambda ).
\end{eqnarray*}%
With results from Maple, it is possible using an incomplete form of
induction to give general expressions for the integrals $I_{\nu
}^{\,(1)}(a,\lambda )$ with even and odd values of $\nu .$ The integrals are
found to contain the Anger functions \cite{Anger} $\mathbf{J}_{1/2}(z)$ and $%
\mathbf{J}_{3/2}(z)$ i.e.%
\begin{equation*}
I_{2n}^{\,(1)}(a,\lambda )=\frac{1}{\sqrt{a}\,2^{\,\epsilon (n)+n}}\left\{ 
\begin{array}{c}
2\sqrt{\frac{2\pi a}{\lambda }}\sum_{k=0}^{2n-1}a_{k,2n}\left( \frac{\lambda 
}{a}\right) ^{k} \\ 
+2\pi \cos (\lambda /a)\sum_{k=0}^{n-1}c_{k,2n}\left( \frac{\lambda }{a}%
\right) ^{2k} \\ 
+4\pi \sin (\lambda /a)\sum_{k=0}^{n-1}d_{k,2n}\left( \frac{\lambda }{a}%
\right) ^{2k+1} \\ 
-\pi \sqrt{\frac{2\pi a}{\lambda }}\mathbf{J}_{1/2}(\lambda
/a)\sum_{k=0}^{n}e_{k,2n}\left( \frac{\lambda }{a}\right) ^{2k} \\ 
+\pi \sqrt{\frac{2\pi a}{\lambda }}\mathbf{J}_{3/2}(\lambda
/a)\sum_{k=0}^{n-1}f_{k,2n}\left( \frac{\lambda }{a}\right) ^{2k+1},%
\end{array}%
\right\}
\end{equation*}%
where the Greubel eta function $\epsilon (n)$ is defined here as 
\begin{equation*}
\epsilon (n)=\lfloor \sqrt{2}n\rfloor +\lfloor \sqrt{3/2}n\rfloor \,\ ,
\end{equation*}%
and \ $\lfloor z\rfloor $ is the floor function. \ The sequence of integers
produced by the eta function has been studied by G. C. Greubel and others %
\cite{Greubel}. The expression for $I_{2n+1}^{\,(1)}(a,\lambda )$ is%
\begin{equation*}
I_{2n+1}^{\,(1)}(a,\lambda )=\frac{1}{\sqrt{a}\,2^{\,\epsilon (n+1)+n}}%
\left\{ 
\begin{array}{c}
2\sqrt{\frac{2\pi a}{\lambda }}\sum_{k=0}^{2n-1}a_{k,2n+1}\left( \frac{%
\lambda }{a}\right) ^{k} \\ 
+2\pi \cos (\lambda /a)\sum_{k=0}^{n}c_{k,2n+1}\left( \frac{\lambda }{a}%
\right) ^{2k} \\ 
+2\pi \sin (\lambda /a)\sum_{k=0}^{n-1}d_{k,2n+1}\left( \frac{\lambda }{a}%
\right) ^{2k+1} \\ 
-\pi \sqrt{\frac{2\pi a}{\lambda }}\mathbf{J}_{1/2}(\lambda
/a)\sum_{k=0}^{n}e_{k,2n+1}\left( \frac{\lambda }{a}\right) ^{2k} \\ 
+\pi \sqrt{\frac{2\pi a}{\lambda }}\mathbf{J}_{3/2}(\lambda
/a)\sum_{k=0}^{n}f_{k,2n+1}\left( \frac{\lambda }{a}\right) ^{2k+1}%
\end{array}%
\right\} .
\end{equation*}%
Alternatively, in keeping with results given above the Anger functions can
be written in terms of the Fresnel functions appearing in $I_{1}^{\,(2)},$
and $\ I_{2}^{\,(2)}$, using the relations%
\begin{eqnarray*}
\mathbf{J}_{1/2}(z) &=&\sqrt{\frac{2}{\pi z}}[C(\sqrt{\frac{2z}{\pi }}%
)\{\sin (z)+\cos (z)\}+S(\sqrt{\frac{2z}{\pi }})\{\sin (z)-\cos (z)\}], \\
\mathbf{J}_{3/2}(z) &=&-\frac{2}{\pi z}+\sqrt{\frac{2}{\pi z}}C(\sqrt{\frac{%
2z}{\pi }})\{\sin (z)\,(1+1/z)-\cos (z)\,(1-1/z)\} \\
&&-\sqrt{\frac{2}{\pi z}}S(\sqrt{\frac{2z}{\pi }})\{\sin (z)\,(1-1/z)+\cos
(z)\,(1+1/z)\}.
\end{eqnarray*}

The coefficients $a_{k,\nu },c_{k,\nu },d_{k,\nu },e_{k,\nu },f_{k,\nu }$
occurring in the $I_{\nu }^{\,(1)}$ integrals are interrelated by (21) from
which we obtain the connection formulas%
\begin{eqnarray*}
\frac{4n}{2^{\Delta \,\epsilon (n)}}a_{2n,2n+1} &=&-e_{n,2n},\hspace{0.25in}%
(k=n), \\
\frac{4n}{2^{\Delta \,\epsilon (n)}}a_{2k,2n+1}+(4k-4n)\,a_{2k,2n}
&=&-e_{k,2n},\hspace{0.25in}(k<n), \\
\frac{4n}{2^{\Delta \,\epsilon (n)}}a_{2k+1,2n+1}+(4k+2-4n)\,a_{2k+1,2n}
&=&f_{k,2n},\hspace{0.25in}(k<n), \\
&& \\
\frac{n}{2^{\Delta \,\epsilon (n)}}c_{n,2n+1} &=&-d_{n-1,2n},\hspace{0.25in}%
(k=n), \\
\frac{4n}{2^{\Delta \,\epsilon (n)}}c_{k,2n+1}+(4k+1-4n)\,c_{k,2n}
&=&-4\,d_{k-1,2n},\hspace{0.25in}(k<n), \\
&& \\
2n\,d_{k,2n+1}+(4k+3-4n)\,d_{k,2n} &=&c_{k,2n,}\hspace{0.25in}(k\leq n-1), \\
&& \\
\frac{4n}{2^{\Delta \,\epsilon (n)}}e_{0,2n+1} &=&(4n-1)e_{0,2n},\hspace{%
0.25in}(k=n), \\
\frac{4n}{2^{\Delta \,\epsilon (n)}}e_{k,2n+1}+(4k+1-4n)\,e_{k,2n}
&=&2f_{k-1,2n},\hspace{0.25in}(k>0), \\
&& \\
\frac{4n}{2^{\Delta \,\epsilon (n)}}f_{n,2n+1} &=&-2e_{n,2n},\hspace{0.25in}%
(k=n), \\
\frac{4n}{2^{\Delta \,\epsilon (n)}}f_{k,2n+1}+(4k-1-4n)\,f_{k,2n}
&=&-2e_{k,2n},\hspace{0.25in}(k<n),
\end{eqnarray*}%
with $\Delta \epsilon (n)=\epsilon (n+1)-\epsilon (n).$ \ These relations do
not appear to be useful except as an internal check on the values of the
coefficients.

\end{document}